\documentstyle{amsppt}

\NoBlackBoxes

\def \prt#1 { {\partial \over \partial
{#1   }}  }
\def \prtf #1 #2{ {\partial #1 \over
\partial {#2   }  }   }

\def \po {p_0} 
\def \cn {\Bbb C^n}

\def \om {\Omega}
\def \z {\zeta}

\def \a {\alpha}
\def \r {\rho}

\def \co {\Cal O}

\def \cu {\Cal U}

\def \J {\text {\rm Jac  }}

\def\-{\overline}

\def \pa {\partial}
\def \hol{\Cal H(\om)}
\def \holp{\Cal H(\om')}

\def \cx{\Bbb C}

\def \b {\beta}

\def \e {\epsilon}

\def \tcpm {T^c_pM}

\def\d {\delta}

\document

\topmatter
\author  M.~S. Baouendi, Xiaojun Huang,
and Linda Preiss Rothschild
\endauthor
\title Nonvanishing of the
differential of holomorphic mappings
at boundary points
\endtitle
\address Department of
Mathematics-0112, University of
California, 
La Jolla, CA 92093-0112
\endaddress
\email sbaouendi\@ucsd.edu
 \endemail
\address  Department of Mathematics,
University of Chicago, Chicago, IL
60637
\endaddress
\email
xhuang\@math.uchicago.edu\endemail
\address Department of
Mathematics-0112, University of
California, 
La Jolla, CA 92093-0112
\endaddress
\email  lrothschild\@ucsd.edu \endemail

\thanks Received April 17, 1995.
\endthanks
\thanks The first and third authors
were partially supported by National
Science Foundation Grant DMS 9203973. 
\endthanks
\leftheadtext {M.~S.~Baouendi, X.~Huang,
and L.~P.~Rothschild}
\rightheadtext{differential of
holomorphic mappings}
\endtopmatter 
\heading  \bf \S 0 Introduction
\endheading

Let $M$ and $M'$ be two smooth
hypersurfaces in $\cn$.  A smooth
mapping
$h: M
\to M'$ is a {\it CR mapping} if its
components are annihilated by the
induced Cauchy-Riemann operator on
$M$. Let $\po
\in M$ and suppose that near $\po$, 
$h$ is the restriction of a holomorphic
mapping
$H$ defined on one side of $M$ near
$\po$ and smooth up to $M$.  We shall
say that
$h$ satisfies the {\it Hopf lemma
property} at
$\po$ if the  component of $H$ normal
to $M_2$ has a nonzero derviative at
$\po$ in the normal direction to $M_1$.
The hypersurface
$M$ is said to be {\it minimal} at $\po
\in M$ if there is no germ of a complex
hypersurface contained in $M$ through
$\po$. Recall the theorem of Tr\'epreau
[T] that if $M$ is minimal at
$\po$, then every CR function defined
on
$M$ near $\po$ extends holomorphically
to at least one side of $M$ in $\cn$
near
$\po$. A stronger condition on a
hypersurface $M$ at a point
$\po$ is that of {\it essential
finiteness}  (as defined in [BJT], 
[BR3], [DA2]).  We will recall this
definition in
\S 1.  We note here that if
$M$ is of D'Angelo finite type at
$\po$ [DA1], then $M$ is essentially
finite at
$\po$ (and hence minimal at $\po$).

  In this paper we prove a general
result of the ``Hopf lemma" type for CR
mappings, with nonidentically vanishing
Jacobians, between real hypersurfaces
in $\cn$. Applications of this result
to finiteness and  holomorphic
extendibility of such mappings are
also given.  The novelty here is that
we make no assumption on the
nonflatness of the mapping or its
Jacobian, nor do we assume that the
hypersurfaces are pseudoconvex or
minimally convex.   

\proclaim {Theorem 1} Let $M$ be a
smooth, connected, orientable
hypersurface in
$\cn$ which is essentially finite at
all points.  Let $h: M \to M'$ be a
smooth CR mapping from $M$ to another
smooth hypersurface $M' \subset \cn$,
with   $\, \J h
\not\equiv 0$. Let $\po \in M$, and
suppose that 
$h^{-1}(h(\po))$ is a compact subset of
$M$.  Then $h$ satisfies the Hopf lemma
property at $\po$.   
\endproclaim

If the hypersurfaces are pseudoconvex,
the result above follows from the
classical Hopf lemma for harmonic
functions, as proved in Fornaess [F]
(see also [BBR]). Other results of the
Hopf lemma type for CR mappings were
previously obtained in [BR3] and
[BR5]. 

As shown in [BBR], [BR1], and [DF], 
the Hopf lemma property can be used to
prove holomorphic extension of CR
mappings between real analytic
hypersurfaces.  From Theorem 1 we
obtain the following corollary.

\proclaim {Corollary 1} Let $M$, $M'$,
$h$ be as in Theorem {\rm 1}, and
assume in addition that
$M$ and $M'$ are real analytic.  Then
$h$ extends holomorphically to a
neighborhood of $\po$ in $\cn$.
\endproclaim

In the global case, i.e. when $M$ and
$M'$ are compact boundaries, the
 compactness of
$h^{-1}(h(\po))$ is automatically
satisfied, yielding the following
result. 
\proclaim {Theorem 2} Let $\om$ and
$\om'$ be bounded domains in
$\cn$ with smooth boundaries, such that
$\partial \om$ is essentially finite at
all points.  Suppose
$H:\om\to
\om'$ is a proper holomorphic mapping,
smooth up to $\partial \om$. Then
$H$ satisfies the Hopf lemma property
at every point $p \in
\partial
\om$. Furthermore,
$H$ is finite-to-one on $\- \om$.  
\endproclaim

When $\om$ and $\om'$ in Theorem 2 are
real analytic, we obtain a new proof of
the following result of the second
author and Pan [HP], extending 
earlier results in [BR1], [DF],
[BR3].

\proclaim {Corollary 2} If
$\om$ and $\om'$ are bounded domains in
$\cn$ with real analytic boundaries,
and 
$H: \om \to \om'$ is a proper
holomorphic mapping, smooth up to
$\partial
\om$,  then $H$ extends
holomorphically to a neighborhood of
$\-
\om$ in $\cn$. 
\endproclaim

It should be noted that Theorem 2 and
Corollary 2 may be proved more directly
(see Remark 2.2).

Another application of Theorem 1 is a
propagation result of the Hopf lemma
property (Theorem 3), as well as real
analyticity (Corollary 3), analogous to
the classical Hartog's theorem for
extension of holomorphic functions.

\proclaim {Theorem 3}  Let $M$ be a
smooth, orientable, connected
hypersurface in
$\cn$ which is essentially finite at
all points, and let $h: M
\to M'$ be a smooth CR mapping from
$M$ to another smooth hypersurface
$M'\subset \cn$, with 
$\J h
\not\equiv 0$. Suppose that
$U_1$ and $U$ are relatively compact
open subsets of $M$, with $\-U_1
\subset U$.  Then if the Hopf lemma
property holds at every point in
$U\backslash U_1$, it also holds
everywhere in $U$.
\endproclaim 

\proclaim {Corollary 3}  Let $M$ be a
real analytic, orientable, connected
hypersurface in
$\cn$ which is essentially finite at
all points, and let $h: M
\to M'$ be a smooth CR mapping from
$M$ to another real analytic
hypersurface
$M'
\subset \cn$, with $\J h
\not\equiv 0$. Suppose that $U_1$ and
$U$ are relatively compact open
subsets of $M$, with $\-U_1
\subset U$. If $h$ is real analytic in
$U\backslash U_1$, then $h$ is real
analytic everywhere in $U$ and hence
extends holomorphically to an open
neighborhood of $U$ in $\cn$.
\endproclaim

\remark {Remark {\rm 0.1}} In Theorem 1
and Corollary 1, the condition $\J h
\not\equiv 0$ may be replaced by the
stronger condition that $M'$ does not
contain any nontrivial complex variety
through $\po$.  (See e.g. [BR4].) A
similar statement holds for Theorem 3
and Corollary 3.
 \endremark

Some of the results of the
present paper, including Corollary 1,
were announced earlier by the
second author. Also,  a recent preprint
of Y. Pan [P] contains a special case 
of Corollary 1 above and other related
results.

\heading  \bf \S 1
Preliminaries\endheading

Let $M$ be a smooth real hypersurface
in
$\cn$.  For $p \in M$, we denote by
$T_pM$ the real tangent space of $M$ at
$p$ and by $\cx T_pM$ its
complexification.  We denote by $\Cal
V_pM$ the complex subspace of
$\cx T_pM$ consisting of all
antiholomorphic vectors tangent to
$M$ at $p$, and by $\tcpm = \text
{{\rm Re }}\Cal V_pM$ the complex
tangent space of $M$ at $p$ considered
as a real subspace of $T_pM$. Note
that if
$h$ is a smooth CR map from
$M$ to a hypersurface  $M'$, then $h$
satisfies the Hopf lemma property
mentioned in \S 0 if and only if
$$dh(T_{\po} M) \not\subset T^c_
{h(\po)} M'.$$  Note that for this
form of the definition, it is not
necessary to assume that $h$ extends
holomorphically to one side of
$M$.

 If
$\rho(z,\-z)$ is a defining function
for
$M$ near $\po =0$, with
$\rho(0) = 0$ and $d\rho (0)\not= 0$,
we consider the formal Taylor series
of $\r$ in $z$ and $\-z$  at $0$ and
write
$R(z,\z)$ for its complexificiation,
i.e.
$R(z,\z) \equiv \sum
a_{\a,\b}z^\a\z^\b$, where $\a! \b!
a_{\a\b}=
\r_{z^\a\-z^\b}(0)$. Let $X_1, \ldots,
X_n$,  be the vector fields with
formal power series coefficients given
by
$$X_j = R_{\z_n}(0,\z){\pa \over \pa
\z_j} - R_{\z_j}(0,\z){\pa \over \pa
\z_n},
\ \ j=1,\ldots, n-1, $$ where we have
assumed $\r_{z_n}(0)\not= 0$.     For a
multi-index
$\alpha = (\a_1,\ldots,\a_{n-1})$ we
define
$c_\a(z)$ in  the ring of convergent
power series in $n$ complex variables,
by 
$$ c_\a(z)= X^{\a} R(z, \z)|_{\z=0},
\tag 1.1
$$ where $X^\a = X_1^{\a_1}\cdots
X_{n-1}^{\a_{n-1}}$. We say that $M$ is
{\it essentially finite} at
$0$ if the ideal
$(c_\a(z))$ generated by the $c_\a(z)$
in the ring $\Bbb C[[z]]$ of formal
power series is of finite codimension.
It should be noted that this
definition is independent of the
choice of coordinates and defining
function $\r$;  it is given in a
slightly different form in [BR3]. 
Note also that if $M$ is essentially
finite at $\po$, then $M$ is minimal at
$\po$, and that if $M$ is of D'Angelo
finite type, then it is essentially
finite.

Recall that  an {\it analytic disc} in
$\cn$ is a continuous mapping $A:
\-\Delta
\rightarrow \cn $ which is holomorphic
in
$\Delta$, where $\Delta
$ is the open unit disc in the plane.
We say that $A$ is  {\it attached} to 
$M$ if  $A(\partial\Delta) \subset M$. 
  Let $M$ be a smooth hypersurface
 minimal at $\po$. As in [BR5], we  say
that
$M$ is {\it minimally convex at
$\po$} if $M$ is minimal at $p_0$, and
there is a neighborhood $U$ of
$\po$ in
$M$ and a side of the hyperplane
$T_{p_0}M$ in $\cn$ such that the real
derivatives ${\partial \over  \partial 
\xi}
   [A(\xi)]_{|_{\xi = 1}}$  lie on that
side or in $T_{p_0}M$, for all
sufficiently smooth analytic discs $A$
attached to
$U$ with  $A(1)= \po$.  Here
$\z = \xi + i \eta$, with $\z \in
\Delta$.

For the convenience of the reader we
begin by stating a number of known
results, Theorems A, B, C and D
below,  which will be important for
the proofs of Theorem 1.  Theorem A is
a consequence of a result of Tumanov
[Tu], as observed in  [BR5] .

\proclaim {Theorem A [Tu, BR5 ]}
 Let $M$ be a smooth, real
hypersurface in
$\cn$, and assume that $M$ is minimal
at
$\po$.  Then one of the following two
conditions holds.
\roster
\item  
$M$ is minimally convex at $\po$. 
\item Every CR function defined in a
neighborhood of $\po$ in
$M$ extends holomorphically to a full
neighborhood of $\po$ in
$\cn$.
\endroster
\endproclaim

\proclaim {Theorem B [BR5]} If $h$ is a
smooth CR mapping from a smooth
hypersurface $M$ to another smooth
hypersurface
$M'$, with
 $M$  minimal at $\po$, $\J  h
\not\equiv 0$, and
 $M'$  minimally convex at
$p'_0=h(p_0)$,  then the Hopf lemma
property holds at
$\po$. 
\endproclaim

We need also to recall a result which
follows  from Theorem 4 in [BR3].
\proclaim {Theorem C [BR3]}
  Let $H$ be a holomorphic map defined
in a neighborhood of a smooth
hypersurface
$M$ essentially finite at
$\po$, $H(M) \subset M'$, with $M'$
another smooth hypersurface of $\cn$,
and
$\J H \not\equiv 0$.  Then H satisfies
the Hopf lemma property at $\po$ and
$H$ is finite-to-one in a neighborhood
of $\po$.
\endproclaim

We shall also need  a stronger version
of this result, and we indicate here
its proof.

\proclaim {Theorem D}
  Let $H$ be a holomorphic map defined
on one side of a smooth hypersurface
$M$ essentially finite at $\po$, with
$H$ smooth up to $M$.  Suppose
$H(M)
\subset M'$, with $M'$ another smooth
hypersurface of
$\cn$. Then
$\J H$ is not flat at $\po$ if and
only if
$H$ satisfies the Hopf lemma property
at
$\po$.  In addition, if either of these
equivalent conditions is satisfied,
then any smooth extension of
$H$ to a sufficiently small
neighborhood of
$\po$ in $\cn$ is finite-to-one.
\endproclaim

\demo {Proof of Theorem D} If $\J H$ is
not flat at $\po$, we conclude that if
$G$ is a formal transversal component
of
$H$ (as defined in [BR3]), then
$G\not\equiv 0$.  Hence, by Theorem 4
of [BR3], it follows that $H$
satisfies the Hopf lemma property and
is of finite multiplicity. Conversely,
if $H$ satisfies the Hopf Lemma
property at
$\po$, by Theorem 4 of [BR3] it follows
again that $H$ is of finite
multiplicity and also that $M'$ is
essentially finite at $H(\po)$.  By
Theorem 3 of [BR3], we conclude also
that
$H$ is not totally degenerate at
$\po$, in the sense of [BR3], and
hence, using again the Hopf lemma
property, $\J H$ is not flat at
$\po$. 

We may assume $\po=H(\po) = 0$.  Since
$H$ is holomorphic on one side of $M$
and smooth up to the boundary, its
Taylor series at $0$ defines a formal
(not necessarily convergent)
holomorphic map
$H=(\sum a^1_\a z^\a, \ldots, \sum
a^n_\a z^\a)$. The equivalent
conditions above imply that 
$H$ is finite as a formal map. That
is, the ideal generated by the $\sum
a^j_\a z^\a, \ j=1,\ldots,n$,  is of
finite codimension in the ring $\Bbb
C[[z]]$ of formal power series in
$z$. Since the Taylor series of
$H$ coincides with that of any smooth
extension of $H$ to $\cn = \Bbb
R^{2n}$, we conclude e.g. by [GG],
[EL], that this extension is
finite-to-one near
$0$.
  $\square$ 
\enddemo

\heading \S 2 Inverse image of a
nonminimally convex point.\endheading

In this section we shall state and
prove a new result, Theorem 4 below,
from which Theorem 1 will follow.

\proclaim {Theorem 4}  Let $h: M \to
M'$ be a CR map, where $M$, $M'$, 
$h$, $\po \in M$ satisfy all the
conditions of Theorem {\rm 1}. In
addition, suppose that $M'$ is not
minimally convex at $\po'=h(\po)$. 
Then all CR functions on $M$ extend
holomorphically to a full neighborhood
of
$\po$ in $\cn$.  In particular, $h$
extends holomorphically to a
neighborhood of
$\po$ and satisfies the Hopf lemma
property at $\po$.  
\endproclaim

Before proving Theorem 4, we note that
Theorem 1 is a consequence of Theorem 4
and Theorem B above.  Indeed, if
$M'$ is  minimally convex at $h(\po)$,
then
 since any essentially finite
hypersurface is minimal at all points,
Theorem 1 follows from Theorem B. On
the other hand, if $h(\po)$ is not
minimally convex, Theorem 1 is an
immediate consequence of Theorem 4.  

In the rest of this section, we shall
prove Theorem 4. We may assume that
$\po = p'_0 = 0$, and we let
$Z_M = h^{-1}(0)$.  Note that $Z_M$ is
a compact subset of $M$ by the
assumptions of the theorem.  Hence
without loss of generality, we shall
assume that
$M$ is bounded.

The following lemma shows that we can
reduce the proof of the theorem to the
case where $h$ extends holomorphically
to one side of $M$.

\proclaim {Lemma 2.1} Under the
assumptions of Theorem {\rm 1}, there
exists an open neighborhood $\cu$ of
$0$ in
$\cn$ such that $Z_M
\cap \cu$ is compact in $M \cap \cu$
and
$h$ extends holomorphically to $\cu^+$,
one side of $M$ in $\cu$.
\endproclaim

\demo {Proof} Since $M$ is essentially
finite and hence minimal at all
points, it follows that $h$ extends
holomorphically to at least one side
of $M$ at each point.  Since
$M$ is orientable, it is given by a
global smooth defining function
$\rho$ with nonvanishing gradient on
$M$. We may assume that 
$h$ extends to the plus side of
$M$, ( i.e. where $\rho (z) > 0$) near
$0$. Let $M_1$ be the largest connected
open subset of $M$ containing $0$ such
that $h$ extends holomorphically to the
plus side of $M$ near every point of
$M_1$.

If $M_1 = M$, then the Lemma is an
immediate consequence of the
assumptions of Theorem 1.  Assume
therefore that
$M_1$ is a proper subset of $M$ and let
$\pa M_1$ be its boundary in $M$.  For
$\d > 0$, let $M_1^\d = \{p \in M_1:
\text {dist}(p,
\partial M_1) > \d\}$. Since at every
point of $M$, $h$ extends
holomorphically to at least one side
of $M$, it follows from the definition
of $M_1$ that there is an open
neighborhood $U$ of
$\partial M_1$ in $M$ such that $h$
extends holomorphically to both sides
of
$M$ at every point in $U \cap M_1$.
Applying Theorem C, we conclude that
$Z_M
\cap U \cap M_1$ is a discrete set.  

Let $\partial M$ be the boundary of
$M$ in
$\cn$ and choose $a > 0$ such that
$a < \text {dist}(Z_M,\pa M)$ (which is
possible by the assumption of the
theorem).  Denote by
$M^{a}  =
\{p
\in M: \text {dist}(p,
\partial M) > a\}$. Note that $\pa M_1
\cap \-{M^a}$ is compact in
$M$.  Therefore, there exists $\d_0 >0$
such that for all $\d$, $0 < \d <
\d_0$, we have
$$
\pa M_1^\d \cap Z_M = \-{M^a} \cap \pa
M_1^\d \cap Z_M \subset U \cap M_1.
$$ By compactness and the discreteness
mentioned above, we conclude that
$\pa M_1^\d
\cap Z_M$ is a finite set.  Since these
sets are all disjoint for different
$\d$'$s$, we conclude that there exists
$\d_1$, with $0 < \d_1 <
\d_0$, such that $\pa M_1^{\d_1} \cap
Z_M =
\emptyset$.  Now the lemma follows by
taking $\cu$ to be a sufficiently small
open neighborhood of
$M_1^{\d_1}$ in $\cn$.
  $\square$ 
\enddemo

By Lemma 2.1, after shrinking $M$
if necessary, we may now assume that
there is a connected open set
$\co$ in
$\cn$ such that:
\roster 
\item "(i)" $\co \cup M$ is a
manifold with boundary of class
$C^\infty$.  
\item "(ii)" $h$ extends
holomorphically to
$\co$; if $H$ denotes the holomorphic
extension of  $h$, then $H
\in C^\infty(\-\co)$.
\item "(iii)" $h(0) = 0$.
\item "(iv)" $Z_M = h^{-1}(0)$
is a compact subset of $M$.
\endroster

We write 
$$Z = H^{-1}(0)\cap \co.$$ We shall
show that we can take $H$ to be a
proper mapping of an open domain to
its image.  The following lemma is
crucial in this construction.

\proclaim {Lemma 2.2} Let $V$ be a
connected open neighborhood of $Z_M$ in
$M$, with $\-V$ a compact subset of
$M$.  For $\d > 0$, let
$$ \co^\d = \{ z \in \co: \text {\rom
{dist}} (z,\-V) < \d\},
\tag 2.3$$  and $\pa \co^\d = S_1^\d
\cup S_2^\d$, with $S_1^\d = \pa
\co^\d \cap M$ and $S_2^\d = \pa\co^\d
\backslash S_1^\d$.  Then for any
$\d_0 > 0$ there exists $\d$,
$0 < \d < \d_0$,  such that
$S_2^\d \cap Z =
\emptyset$.\endproclaim 

\demo {Proof} Note that by assumption,
$\-{S_2^\d} \cap Z_M = \emptyset$. 
Hence there exists $\e > 0$ such that
for any
$\d$ sufficiently small,  
$$ Z \cap \{z \in S_2^\d:\  \text {\rom
dist} (z,\-{S_2^\d}\cap M) < \e\} =
\emptyset. \tag 2.4 $$ Let $Z' = \{ z
\in Z: z\ \text {is not an isolated
point of}
\  Z\}$. If there exists $\d > 0$ such
that
$Z' \cap
\co^\d =
\emptyset$, then $Z \cap \co^\d$ is
countable, and the conclusion of the
lemma follows since the sets $S^\d_2$,
as
$\d$ varies, are disjoint. 

To complete the proof of the lemma, we
shall assume
$$Z' \cap \co^\d \not= \emptyset, \tag
2.5$$
 for all $\d$ sufficiently small, and
reach a contradiction.  It is clear
that under condition (2.5) we have
$$\overline { Z'} \cap M \not=
\emptyset.
\tag 2.6$$  We claim that $\J H$
vanishes to infinite order at every
point
$p \in  \overline { Z'} \cap M$.
Indeed if
$\J H$ does not vanish to infinite
order at a point $p \in Z_M$, then by
Theorem D,
$p$ is an isolated point of
$H^{-1}(0)$ in
$\-\co$.  Since this cannot be the case
for $p \in \overline { Z'} \cap
M\subset Z_M$, the claim follows. 

Now let $T_2^\d = \{z \in S_2^\d:\
\text {\rom {dist}} (z,\overline
{S_2^\d}\cap M)
\ge
\e\}$, where $\e$ satisfies (2.4). 
Using (2.4), we note that
$$\overline { Z'}\cap \pa \co^\d
\subset T_2^\d \cup (Z_M
\cap
\overline {Z'}). 
$$ Since $T_2^\d$ is compactly
contained in $\co$ for sufficiently
small $\d$, and
$\J  H$ is holomorphic in
$\co$, there exists $C > 0$ such that
for all multi-indices $\a$
$$ \sup_{z \in T_2^\d} |D^\a \J H(z)|
\le C^{|\a|+1}\a!. \tag 2.7
$$   By the maximum principle on
complex varieties (see e.g. [N1]) we
have, 
$$\sup_{z \in \overline {Z'}\cap
\-{\co^\d}} |D^\a
\J H(z)| =
\sup_{z \in \overline { Z'}\cap \pa
\co^\d} |D^\a \J H(z)|. \tag 2.8$$
However, as proved in the claim above,
$\J H$ vanishes to infinite order on
$\overline { Z'} \cap M$.  Hence, in
view of (2.7) and (2.8)
$$\sup_{z \in \overline {Z'}\cap
\-{\co^\d}} |D^\a
\J H(z)| =
\sup_{z \in T_2^\d} |D^\a \J H(z)|\le
C^{|\a|+1}\a!. \tag 2.9$$ This proves
that the radius of convergence of $\J
H (z), z
\in Z'$, is greater than a positive
constant which is independent of the
distance to
$M$.  Hence
$\J H$ extends holomorphically to a
full neighborhood in $\cn$ of each
point of 
$\overline { Z'}
\cap M$. Since
$\J H$ vanishes to infinite order
there, it follows that $\J H \equiv
0$, contrary to assumption. We
conclude  that (2.6), and hence (2.5),
cannot hold, which completes the proof
of Lemma 2.2.
$\square$
\enddemo

In reducing the proof of Theorem 4 to
the global case of a proper mapping we
shall use the following. 

\proclaim {Proposition 2.10} Let $M$
be a
 connected  hypersurface of class $C^0$
with
$M
\subset \pa \co$, where $\co$ is an
open bounded domain in $\cn$, and let
$H$
 be a holomorphic mapping in
$\co$, continuous up to the boundary,
$\J H
\not\equiv 0$,  with
$H(M)$ contained in another
hypersurface
$M'$ of class $C^0$.  Suppose $0 \in
M$,
$H(0) = 0$, and 
$$H^{-1}(0) \cap \pa \co  \subset M. 
\tag 2.11$$ Then there is a subdomain
$\co_1
\subset \co$ satisfying
\roster
\item "(i)" $0 \in \pa \co_1$, and
there exists a sequence $\{z_j\}
\subset 
\co_1$  such that $z_j \to 0$ and
$H(z_j)$ stays strictly on one side of
$M'$;  
\item "(ii)" there exists $U$, a
neighborhood of
$0$ in $M'$, with $\overline
{H({\co_1})}\supset U$;
\item "(iii)" $H: \co_1\to H(\co_1)$
is a proper map.
\endroster
\endproclaim

\demo {Proof} We begin with the
following lemma, which describes a
well-known construction, see e.g. [BC].

\proclaim {Lemma 2.12} Let $\co \subset
\cn$ be an open bounded domain, and
suppose $H: \co \to  \cn$ is a
holomorphic mapping, continuous up to
$\pa \co$.  Let  
$$D = \{ z \in \co: H(z) \not\in H(\pa
\co)\}. $$ If $D \not= \emptyset$, then
$H: D \to H(D)$ is finite-to-one and
hence open.  Furthermore, if
$D'_1$ is any connected component of
$H(D)$, and $D_1$ a connected
component of
$H^{-1}(D_1')$, then $H: D_1 \to D_1'$
is a proper map.
\endproclaim
\demo {Proof} Since this result is in
the ``folklore", we shall be brief. 
We  assume $D \not= \emptyset$.  If
$H$ is not finite-to-one, there exists
$w \in H(D)$ for which $H^{-1}(w)$ has
an accumulation point $z_0$ in $\-D$
(and hence $H(z_0) = w$).  By the
definition of $D$,
$$ H(D) \cap H(\pa \co) = \emptyset, \
\ 
\text {and} \ \ \ H(\pa D)
\subset H(\pa \co). \tag 2.13$$  Hence
$z_0
\not\in \pa D$.  On the other hand, if
$z_0 \in D$, then there is a nontrivial
variety contained in $H^{-1}(w)$, which
would necessarily intersect $\pa D$. 
Since this is also impossible, by the
definition of $D$,
$H$ is finite-to-one and hence open
(see e.g. [R]).
 
 To show that $H: D_1 \to D_1'$ is
proper, suppose $z_j \to z_0$, $z_j
\in D_1$, $z_0 \in \pa D_1$.  Then by
continuity
$H(z_j)
\to H(z_0)=w_0 \in
\-{D_1'}$. We claim that $w_0 \in \pa
D_1'$. Indeed, if $w_0$ is an interior
point of $D_1'$, let $V'$ be an
neighborhood of $w_0$ in $D_1'$.
Consider
$H$ as a map from
$\-\co$ to $\cn$.  Then a component of
$H^{-1}(V')$ is contained in $D_1$, by
the definition of $D$.  Then
$z_0$ would be an interior point of
$D_1$, contrary to assumption.  This
proves Lemma 2.12.   
  $\square$ 
\enddemo 

We may now complete the proof of
Proposition 2.10. Let $\rho' $ be a
defining function for $M'$  near $0$. 
Without loss of generality, we may
assume that there exist $z_j \in \co,
j=1,2,\ldots,$ with 
$$ \lim_j z_j = 0 \ \ \text {and} \ \
\rho'(H(z_j)) >0.
\tag 2.14 $$ Indeed, we first select
$z_j
\in \co$ with $\J H(z_j) \not= 0$. 
Since
$H$ is open near such a $z_j$, by
slightly moving $z_j$ if necessary, we
may assume
$H(z_j) \not\in M'$.  Replacing $\rho'$
by 
$-\rho'$ and selecting a subsequence if
necessary, we reach the desired
conclusion (2.14).

 Let $D$ be as in Lemma 2.12. By
hypothesis (2.11) and the continuity of
$H$ it follows that
$H(\pa
\co
\backslash M)$ is a compact set which
does not contain $0$. Hence, by taking
$z_j$ sufficiently close to $0$, we
may assume that the points
$z_j$ chosen in (2.14) are in $D$.  We
shall show that there exists $\e > 0$,
arbitrarily small, such that 
$$  \{w\in \cn: |w| < \e \ \text {and}
\
\rho'(w) > 0\} \equiv W_\e \subset
H(D).
\tag 2.15$$
 Suppose that (2.15) is proved. Let
$D'_1$ be the connected component of
$H(D)$ containing the connected open
set
$W_\e$. We claim that there is a
connected component $D_1$ of
$H^{-1}(D'_1)$ such that $0
\in \pa D_1$.  Indeed, by Lemma 2.12,
the restriction of $H$ to $D$ is
finite-to-one, and the restriction to
any connected component of   
$H^{-1}(D'_1)$ is proper, and hence
onto
$D_1'$.  Therefore,
$H^{-1}(D'_1)$ consists of finitely
many connected components
$D_k$.  Choose one of these components,
say $D_1$, which contains infinitely
many of the $z_j$.  Then $0 \in \pa
D_1$.  Since, by Lemma 2.12, the
restriction of
$H$ to
$D_1$ is proper onto
$H(D_1)=D_1'$, Proposition 2.10 will
follow by taking $\co_1 = D_1$.

It remains to prove (2.15). Choose $\e$
such that $$0<\e <\text{dist} (0,
H(\pa\co\backslash M)), \tag 2.16$$ and
such that the open set
$W_\e$ defined in (2.15) is connected.
Let
$j_0$ be such that $H(z_{j_0}) \in
W_\e$.  Let $w \in W_\e$ be arbitrary,
and
$\gamma(t),  0\le t \le 1$, be a
continuous curve connecting
$H(z_{j_0})$ and $w$ and contained in
$W_\e$.  Assume by contradiction that
$w \not\in H(D)$.  Since
$H(D)$ is open, there exists $t', 0 <
t'
\le 1$, such that $\gamma(t)
\in H(D)$ for $0 \le t < t'$, but
$\gamma(t') \not\in H(D)$.  Now choose
a sequence
$t_k < t'$, with $t_k \to t'$, and $p_k
\in D$ with $H(p_k) =
\gamma(t_k)$ and $p_k \to p' \in \-D$. 
Since $H(p') = \gamma(t')
\not\in H(D)$, it follows from the
definition of $D$ that $p' \in \pa D$. 
Recall that $H(\pa D) \subset H(\pa
\co)$.  Hence $H(p') \in H(\pa \co)$.
In view of (2.16) and the fact that
$H$ maps
$M$ into $M'$, we must have
$H(p') \in M'$. We reach a
contradiction, since
$H(p')=\gamma(t') \in W_\e$. The proof
of Proposition 2.10 is now complete.
$\square$  
\enddemo

Let $\om$ be a bounded domain in $\cn$
and
$\po \in \cn$.  Recall that $\po$ is in
the {\it holomorphic hull} of
$\om$ if there is a compact subset $K
\subset \om$ such that $\po \in\-{\hat
K}$, where 
$$  \hat K = \{ z \in \om: |f(z)| \le
\sup_{w \in K} |f(w)|\ \hbox {for all}\
f \in \Cal H (\om)\}. 
$$
Here $\Cal H (\om)$ denotes the
space of all holomorphic functions
in $\om$. 
 We observe that when $\po$ is a
boundary point of
$\om$, then $\po$ is in the holomorphic
hull of $\om$ if and only if any
function in $\Cal H (\om)$ extends
holomorphically to some larger domain
which contains $\po$ as an interior
point.
 
\proclaim {Proposition 2.17} Let $\om$
and
$\om'$ be two bounded domains in $\cn$
and
$H$ a proper holomorphic mapping from
$\om$ to $\om'$.  Suppose that $p_0$
and
$p'_0$ are boundary points of $\om$ and
$\om'$ respectively, and that there is
a sequence $\{z_j\}_{j=1}^\infty
\subset \om$ converging to $p_0$ such
that $\lim_j H(z_j) = p'_0$. Suppose
that any  function in $\holp$ is
bounded on the sequence
$\{H(z_j)\}_{j=1}^\infty$.  Then $\po$
is in the holomorphic hull of $\om$.
\endproclaim

\remark {Remark {\rm 2.18}} Note that
the hypothesis of the proposition is 
satisfied if $\po'$ is in the envelope
of holomorphy of $\om'$.
\endremark
\demo {Proof} By using a standard
result (see e.g. [N2] Chapter 7, Lemma
2), it suffices to prove the following
claim:

{\sl Each function in $\hol$
is bounded on
$\{z_j\}$.  More precisely, for any
$f\in \hol$, there exists a constant
$C_f > 0$ such that $|f(z_j)| \le C_f$
for all $j$.}

To prove the claim, we note that $H$ is
finite-to-one on $\om$ since it is
proper.  Hence there exists
$m$ such that each $w \in \om'$ has $m$
pre-images, $g_k(w), k = 1
\ldots m$, counted with multiplicity
(see e.g. [R]).  Now let $f\in\hol$
and denote by $\sigma_1(w),
\ldots ,\sigma_m (w)$ the elementary
symmetric functions of
$f(g_k(w)),  k = 1, \ldots, m$.  By
well known results (see e.g. [R]) the
$\sigma_k(w)$ are holomorphic in $\om'$
and hence, by hypothesis, uniformly
bounded on the sequence
$\{H(z_j)\}$.  If we let
$w_j = H(z_j)$,
 we observe that 
$f(z_j)$ is one of the roots of the
polynomial $X^m -
\sigma_1(w_j)X^{m-1} +
\ldots + (-1)^m \sigma_m (w_j)$. 
Since the coefficients of this
polynomial are bounded, independently
of $j$, it follows that
 the $f(z_j)$ are bounded,
independently of
$j$. This proves the claim and hence
Proposition 2.17.
  $\square$ 
\enddemo

\demo {Proof of Theorem \text {\rom 4}}
First, we prove that under the
assumptions of Theorem 4, 
$M'$ is minimal at
$\po' =h(p_0)$.  Indeed, suppose not.
Then there is a complex hypersurface
$\Sigma$ contained in
$M'$ through $\po'$.  Hence, this
hypersurface must contain all small
analytic disks $A'$ attached to $M'$
with
$A'(1) =
\po'$. On the other hand since
$M$ is minimal at
$\po$, the boundaries of 
small analytic discs
$A$ attached to
$M$ with $A(1) = \po$ cover a full
neighborhood of $\po$ in $M$ [Tu]. 
Since we can take 
$A' = h\circ A$, this contradicts the
assumption that $\J h \not\equiv 0$.
(See also [E] for related results.)

 By Lemma 2.1, we may assume that $h$
admits a  holomorphic extension  $H$
to  one side of
$M$, and that conditions (i)--(iv)
preceding Lemma 2.2 are satisfied, so
that we may apply Lemma 2.2.   If
$\d$ satisfies the conclusion of Lemma
2.2,  then
$H$ satisfies the hypotheses of
Proposition 2.10     with
$\co =
\co^\d$. We then obtain from
Proposition 2.10 a subdomain $\Cal O_1$ of
$\co$ such that the restriction of $H$
to
$\Cal O_1$ is a proper mapping from $\Cal O_1$ to
$\Cal O_1'$, continuous up to the boundary,
with $\po \in \pa \Cal O_1$ and $\po' =
H(\po)
\in \pa \Cal O_1'$. Moreover, there exists a
sequence $\{z_j\} \subset \Cal O_1$ such
that
$w_j = H(z_j) \to \po'$, with $\{w_j\}$
strictly on one side of
$M'$.

Since $M'$ is  minimal,
but not minimally convex at $\po'$, by
assumption, it follows from Theorem A
that any CR function defined near
$\po' \in M'$ extends holomorphically
to a full neighborhood of $\po'$ in
$\cn$.  Now, by using the Baire
Category Theorem (see e.g. [BR2,
Theorems 7 and 8] for a more general
result) we conclude that there is a
connected neighborhood  $\Cal U'$ of
$\po'$ in $\cn$ with $\Cal U'\cap \Cal O'_1
\not= \emptyset$ such that every
function in
$ \Cal H (\Cal O'_1)$ extends
holomorphically to
$\Cal U'$.  In particular, we see that
any such function is uniformly bounded
on
$\{w_j\}$.   Using
Proposition 2.17 we conclude that
$\po$ is in the holomorphic hull of
$\Cal O_1$, which lies on the side of
$M$ to which every CR function near
$\po$ extends. It follows immediately
that every CR function near $\po$ on
$M$ extends holomorphically to a full
neighborhood of
$p_0$ in $\cn$.  The Hopf lemma
property then follows from Theorem C
above. The proof of Theorem 4 (and
hence that of Theorem 1) is now
complete.   
  $\square$ 
\enddemo
\comment
We remark that an immediate corollary
of Proposition 2.17 is the following,
which may already be known in the
literature.

\proclaim {Proposition 2.19}  Let $\om$
and $\om'$ be two bounded domains in
$\cn$ and let $H$ be a proper
holomorphic map from $\om$ to $\om'$. 
Then $\om$ is a domain of holomorphy
if and only if $\om'$ is a domain of
holomorphy.
\endproclaim
\endcomment 

\heading  \bf \S 3 Consequences of
Theorem 1 and remarks
\endheading

In this section we prove the other
results stated in the introduction and
make some remarks.

We first note that Corollary 1 follows
easily from Theorem 1 and the following
holomorphic extendibility result, which
is a consequence of Theorem 1 of [BR1]:
\proclaim {Theorem E [BR1]} Let
$h: M\to M'$ be a smooth CR map, with
$M$ and $M'$ real analytic
hypersurfaces in
$\cn$.  Assume that
$M$ is essentially finite at
$\po$ and that
$h$ satisfies the Hopf Lemma property
at
$\po$. Then $h$ extends holomorphically
to a full neighborhood of $\po$ in
$\cn$.   
\endproclaim

\demo {Proof of Theorem {\rm 2}} In
order to apply Theorem 1, we note
first that since
$H$ is proper, $\J H \not\equiv 0$ in
$\om$.  Hence its boundary value on
$\pa
\om$ does not vanish identically.  Note
also that  for any
$\po
\in
\pa
\om$,
$H^{-1}(H(\po))$ is closed in $\pa
\om$ and hence compact. We may now
conclude by Theorem 1 that the Hopf
lemma property holds at each point in
$\pa \om$. 

 To prove that $H$ is finite-to-one in
$\-\om$, we observe first that
$H$ is finite-to-one in $\om$, since it
is proper (see e.g. [R]). Since the
Hopf lemma property holds at
$\po$, we may apply the last part of
Theorem D to conclude that  for any
$\po \in \pa
\om$, $H$ is finite-to-one in a
neighborhood of $\po$ in $\-\om$.  The
desired result then follows by
compactness of $\-\om$.  
  $\square$ 
\enddemo

\remark {Remark {\rm 2.1}} It also
follows from Theorem 4 in [BR3] that
under the hypotheses of Theorem 2, $\pa
\om'$ is also essentially finite at all
points. 
\endremark 

\demo {Proof of Corollary {\rm 2}}  By
a result of Diederich and Fornaess
[DF], any compact real analytic
boundary in $\cn$ does not contain a
nontrivial complex variety and hence
is essentially finite. We may then
apply Theorem 2 to conclude that the
Hopf lemma property is satisfied at
every point of $M$. The conclusion of
Corollary 2 then follows from Theorem
E.
  $\square$ 
\enddemo 

\remark {Remark {\rm 2.2}}  In fact,
Corollary 2 may be proved much more
directly by using Proposition 2.17 
together with Theorems A, B, C, and E.
\endremark

\demo {Proof of Theorem {\rm 3}}
 It suffices to show that if $p_0 \in
U_1$, then $h$ satisfies the Hopf lemma
property at $p_0$. By taking the
connected components of $U$ and
$U_1$ containing
$\po$, we may assume, without loss of
generality, that
$U_1$ and
$U$ are connected.  Let $E=\{p \in U:
h(p) = h(\po)\}$.  Since $M$ is
essentially finite, $h$ extends
holomorphically to one side of
$M$ near any point. Therefore, since by
assumption the Hopf lemma property
holds in
$U\backslash U_1$, it follows from
Theorem D that
$E\cap (U\backslash U_1)$ is a discrete
set.  

For $\d > 0$, sufficiently small, let
$$U^\d = \{p \in U: \text{\rm dist}(p,
\pa U) > \d\},$$ and let
$\pa U^\d$ be its boundary.  By the
discreteness established above and the
compactness of $\pa U^\d$, we conclude
that for sufficiently small
$\d$, that
$\pa U^\d \cap E$ is finite.  Hence
there exists $\d_1>0$ for which the set
$\pa U^{\d_1} \cap E$ is empty.  It is
now easy to check that the hypotheses
of Theorem 1 are satisfied for $h$ and
$\po$ by  taking
$M=U^{\d_1}$. This proves Theorem 3.
$\square$ 
\enddemo   

\demo {Proof of Corollary {\rm 3}}
Since
$h$ is real analytic at all points of
$U\backslash U_1$ and $M$ is real
analytic, $h$ extends holomorphically
to a full neighborhood in $\cn$ of
each such point.  By Theorem C, $h$
then satisfies the Hopf lemma property
in all of
$U\backslash U_1$ and hence in all of
$U$ by Theorem 3.  Applying Corollary
1, we then have that
$h$ extends holomorphically to a full
neighborhood of $U$ in $\cn$.
  $\square$ 
\enddemo

\enddocument